\documentclass[letterpaper, 10 pt, conference]{ieeeconf}

\IEEEoverridecommandlockouts
\usepackage{amsmath} 
\usepackage{amssymb}  
\usepackage[utf8]{inputenc}
\usepackage{color}
\usepackage{tikz}
\usetikzlibrary{shapes,arrows.meta}
\usepackage{xspace}
\usepackage{siunitx}
\usepackage{cite}
\usepackage{bm,fixmath}

\usepackage{graphicx}
\graphicspath{{./image/}}

\usepackage{algorithm}
\usepackage{algpseudocode}

\algrenewcommand\algorithmicrequire{\textbf{Input:}}

\DeclareMathOperator*{\argmin}{arg\,min}
\DeclareMathOperator{\prox}{prox}
\DeclareMathOperator{\refl}{refl}

\newtheorem{remark}{Remark}
\newtheorem{proposition}{Proposition}

\newtheorem{assumption}{Assumption}
\newtheorem{lemma}{Lemma}
\newtheorem{corollary}{Corollary}

\newcommand{\norm}[1]{\left\lVert#1\right\rVert}

\newcommand{\reals}{\mathbb{R}}

\newcommand{\uv}{\mathbold{u}}
\newcommand{\x}{\mathbold{x}}
\newcommand{\y}{\mathbold{y}}
\newcommand{\w}{\mathbold{w}}
\newcommand{\z}{\mathbold{z}}

\newcommand{\Am}{\mathbold{A}}

\renewcommand{\Im}{\mathbold{I}}
\newcommand{\Pm}{\mathbold{P}}
\newcommand{\Wm}{\mathbold{W}}

\newcommand{\one}{\pmb{1}}
\newcommand{\zero}{\pmb{0}}
\newcommand{\Ts}{T_\mathrm{s}}
\newcommand{\Nc}{N_\mathrm{C}}
\newcommand{\Np}{N_\mathrm{P}}

\renewcommand{\baselinestretch}{0.96}

\newif\ifarxiv
\arxivtrue 

\title{\LARGE \bf
Distributed Prediction-Correction ADMM for \\ Time-Varying Convex Optimization
}
\author{Nicola Bastianello, Andrea Simonetto, Ruggero Carli
\thanks{N. Bastianello and R. Carli are with the Department of Information Engineering, University of Padova, Italy.\newline {\tt\scriptsize nicola.bastianello.3@phd.unipd.it, carlirug@dei.unipd.it}. \newline A. Simonetto is with IBM Research Ireland. Dublin, Ireland. \newline{\tt\scriptsize andrea.simonetto@ibm.com}.}
}

\begin{document}

\maketitle

\begin{abstract}
This paper introduces a dual-regularized ADMM approach to distributed, time-varying optimization. The proposed algorithm is designed in a prediction-correction framework, in which the computing nodes predict the future local costs based on past observations, and exploit this information to solve the time-varying problem more effectively. In order to guarantee linear convergence of the algorithm, a regularization is applied to the dual, yielding a dual-regularized ADMM. We analyze the convergence properties of the time-varying algorithm, as well as the regularization error of the dual-regularized ADMM. Numerical results show that in time-varying settings, despite the regularization error, the performance of the dual-regularized ADMM can outperform inexact gradient-based methods, \emph{as well as exact dual decomposition techniques}, in terms of asymptotical error and consensus constraint violation.
\end{abstract}

\section{Introduction}\label{sec:intro}
In this paper, we are interested in solving the time-varying optimization problem
\begin{equation}\label{eq:continuous-problem}
	\bar{\x}^*(t) = \argmin_{\bar{\x} \in \reals^n} \sum_{i=1}^N f_i(\bar{\x}; t)
\end{equation}
over a network of $N$ computing and communicating nodes, each one privately storing a term $f_i$ of the cost function. We assume that each of the cost functions $f_i: \reals^n \times \reals_{+} \to \reals$ is a strongly convex and smooth function uniformly in time $t \geq 0$, so that the solution trajectory $\bar{\x}^*(t)$ exists and it is unique.
Problems of the form~\eqref{eq:continuous-problem} have recently attracted an increasing amount of attention, see~\cite{Rahili2015, akbari_distributed_2017, Shahrampour2018, Ye2015, Rahili2015a, Gong2016, Xi2016a, Sun2017, Maros2019, Maros2017, Ling2013, Jakubiec2013, Zavlanos2013, Paper2, bedi_asynchronous_2019, yi_distributed_2020,zhang_online_2020}, and they naturally appear whenever a group of computing and locally communicating entities need to reach a consensus in a cooperative fashion, without revealing private information on their time-varying costs. In time-varying settings, this can be the case in robotics~\cite{Zavlanos2013}, smart grids~\cite{DallAnese2016}, or transportation networks~\cite{Eser2018}. 

What is challenging about~\eqref{eq:continuous-problem} is that the computing nodes: \emph{(i)} know their private functions only up to the current time $t$, while future functions are unknown and need to be \emph{predicted}; \emph{(ii)} they have limited computation and communication capabilities, so that they cannot solve~\eqref{eq:continuous-problem} exactly at each time $t$. These two challenges can be tackled in the framework of time-varying algorithms, where one sets up online algorithms of limited computation that eventually find and track the solution trajectory as time evolves. 

In this paper, we focus on discrete-time algorithms of the prediction-correction kind~\cite{Simonetto2017a, bastianello_primal_2020}, instead of continuous-time ones~\cite{Fazlyab2018, paternain2019prediction}. Discrete-time algorithms sample the problem~\eqref{eq:continuous-problem} at fixed intervals $t_k$, $k \in \mathbb{N}$ with $t_{k+1} - t_k = \Ts$, which yields the following sequence of time-invariant problems 
\begin{equation}\label{eq:sampled-problem}
	\bar{\x}^*(t_k) = \argmin_{\bar{\x} \in \reals^n} \sum_{i=1}^N f_i(\bar{\x}; t_k).
\end{equation}
The idea is then to devise an online algorithm to approximately solve~\eqref{eq:sampled-problem} within the sampling period, and eventually converge to the optimizer trajectory. Specifically prediction-correction algorithms predict how the cost function changes in time and then correct for errors when a new function is acquired at time $t_{k+1}$, see~\cite{Simonetto2017a, bastianello_primal_2020} and reference therein.  Here, we will devise algorithms that tracks the optimal solution trajectory $\{ \bar{\x}^*(t_k) \}_{k \in \mathbb{N}}$ up to a bounded error and that can be deployed in a distributed fashion.

The key novelty of the paper is a new dual-regularized alternating direction method of multipliers (ADMM), which can also be applied to static problems and it is therefore of independent interest. This dual-regularized ADMM extends the line of research on dual regularizations started in the static setting~\cite{Nedic2011,Devolder2011} and continued in time-varying scenarios~\cite{Simonetto2014d,DallAnese2016}. In general, regularizations change the original problems but improve the convergence properties (e.g., the rate) to the regularized optimizer. Whenever the introduced approximation is acceptable with respect to the added benefit (a faster obtained approximate solution, rather than a slower obtained exact one), then regularized algorithms are preferred to exact ones. We show here that in time-varying scenarios one might ``have the cake and eat it too'': since discrete-time time-varying algorithms never deliver exact solutions and the radius of the bounded error at which one converges depends on the convergence rate, then regularizations can both increase the rate and reduce the asymptotical error, if properly designed. This conclusion is in line with current research on algorithm hierarchies in time-varying optimization and their differences with respect to static optimization~\cite{dallanese_optimization_2019,yuan2020can}.  

The contributions of this paper are as follows:

$\bullet$ We develop distributed prediction-correction algorithms that can be deployed on a network of computing and communicating nodes and prove their convergence. The algorithms are based on a novel dual-regularized alternating direction method of multipliers (ADMM). These algorithms extend the ones available in the literature, \emph{e.g.},~\cite{simonetto_dual_2019, bastianello_primal_2020}, since we are the first to employ the ADMM machinery in a distributed prediction-correction setting. 

$\bullet$ We analyze the novel dual-regularized ADMM (which is also of independent interest) theoretically and practically. In theory, we bound the distance of the dual solution of the regularized problem with the non-regularized one and show linear convergence of the former. In practice, we demonstrate how adding a regularization term is beneficial in time-varying settings, in a numerical example. 

\vskip2mm
\footnotesize
\noindent{\bf Notation.} Vectors and matrices are indicated with $\x \in \reals^n$, and $\Am \in \reals^{n \times m}$, respectively. We denote by $\mathcal{G} = (\mathcal{V}, \mathcal{E})$ the undirected, connected graph describing the distributed system. We denote by $\mathcal{N}_i$ the neighborhood of node $i$, and by $d_i := |\mathcal{N}_i|$ its degree. With $d_\mathrm{M}$ we denote the maximum degree in the network.
The Euclidean norm is denoted by $\norm{\cdot}$, the Kronecker product by $\otimes$. The identity matrix is denoted by $\Im$, and $\one$, $\zero$ denote the column vectors of all ones and zeros, respectively. With $\lambda_\mathrm{m}(\Am)$ and $\lambda_\mathrm{M}(\Am)$ we denote the smallest and largest eigenvalues of a matrix $\Am \in \reals^{n \times m}$.
The convex conjugate of a convex, closed and proper function is defined as $f^\star(\w) = \max_\x \left\{ \langle \w, \x \rangle - f(\x) \right\}$.
The indicator function of a non-empty, closed, convex set $\mathbb{X}$ is denoted by $\iota_{\mathbb{X}}(\x)$, with $\iota_{\mathbb{X}}(\x) = 0$ if $\x \in \mathbb{X}$, and $\iota_{\mathbb{X}}(\x) = +\infty$ otherwise.
Given a convex, closed and proper function $f$, we define its proximal operator as $\prox_{\rho f}(\x) = \argmin_\y \{ f(\y) + \norm{\y - \x}^2 / (2\rho) \}$, $\rho > 0$, and the corresponding reflective operator as $\refl_{\rho f}(\x) = 2 \prox_{\rho f}(\x) - \x$.
A function $f: \reals^n \to \reals$ is $m$-strongly convex, for a constant $m \in \reals_+$, iff $f(\x) - \frac{m}{2}\|\x\|^2$ is convex. The function $f$ is said to be $L$-\textit{smooth} if its gradient is $L$-Lipschitz continuous, or equivalently $f(\x) - \frac{L}{2}\|\x\|^2$ is concave. We denote the class of $m$-strongly convex and $L$-smooth functions with $\mathcal{S}_{m,L}(\reals^n)$.

\normalsize

\section{An ADMM Reformulation}\label{sec:ref}

In this section, we reformulate~\eqref{eq:sampled-problem} as a consensus problem, which we will tackle with a new, regularized, time-varying version of ADMM.

\subsection{Consensus problem formulation}
We introduce the local copies (each for every nodes) $\x_i \in \reals^n$, $i = 1, \ldots, N$, of the unknown variable $\bar{\x}$, and equivalently rewrite~\eqref{eq:sampled-problem} as
\begin{subequations}\label{eq:consensus-problem-1}
\begin{align}
	\x^*(t_k) = & \argmin_{\x \in \reals^{nN}} \sum_{i=1}^N f_i(\x_i;t_k) =: f(\x; t_k)\\
	&\text{s.t.} \ \x_i = \x_j \ \text{if} \ (i,j) \in \mathcal{E} \label{eq:consensus-constraints}
\end{align}
\end{subequations}
where $\x^*(t_k) = \x_k^* = [(\x_{1,k}^*)^\top, \cdots, (\x_{N,k}^*)^\top]^\top$ and $\x_{1,k}^* = \ldots = \x_{N,k}^*$ due to the consensus constraints~\eqref{eq:consensus-constraints}. Further introducing the two \emph{bridge variables} $\y_{ij}\in\reals^{n}$ and $\y_{ji}\in\reals^{n}$ for each edge $(i,j) \in \mathcal{E}$, we can rewrite the consensus constraints as
$$
	\x_i = \y_{ij}, \quad \x_j = \y_{ji} \quad \text{and} \quad \y_{ij} = \y_{ji}.
$$

We define now the matrix $\Am \in \reals^{2n|\mathcal{E}| \times nN}$ as
$$
	\Am = \begin{bmatrix}
		\one_{d_1} & \zero_{d_1} & \cdots & \zero_{d_1} \\
		& \ddots & & \\
		\zero_{d_N} & \cdots & \zero_{d_N} & \one_{d_N}
	\end{bmatrix} \otimes \Im_n
$$
and the permutation matrix $\Pm \in \reals^{2n|\mathcal{E}| \times 2n|\mathcal{E}|}$ which swaps $\y_{ij}$ with $\y_{ji}$ (see~\cite{bastianello_asynchronous_2019} for details). Notice that for each edge, there are \emph{two} consensus constraints of the form $\x_i = \y_{ij}$, and so the bridge variables are $2|\mathcal{E}|$. We can then rewrite the consensus constraints as
$
	\Am \x - \y = \zero \quad \text{and} \quad \y = \Pm \y,
$
and thus the consensus problem~\eqref{eq:consensus-problem-1} at time $t_k$ is equivalent to
\begin{subequations}\label{eq:consensus-problem}
\begin{align}
	\x^*(t_k), \y^*(t_k) = & \argmin_{\x \in \reals^{nN}} \sum_{i=1}^N f_i(\x_i;t_k) + \iota_{\ker(\Im - \Pm)}(\y) \\
	&\text{s.t.} \ \Am \x - \y = 0 \label{eq:lin-constraint}
\end{align}
\end{subequations}
where $\ker(\mathbold{M})$ denotes the null-space of matrix $\mathbold{M}$.

It has been shown in \emph{e.g.}~\cite{shi_linear_2014,bastianello_asynchronous_2019}, how problem~\eqref{eq:consensus-problem} for a fixed time $t_k$ can be solved in a distributed fashion, by allowing each node to communicate only with its neighbors and by employing ADMM. Here we look at time-variant versions of ADMM, where only a limited number of steps are allowed at each time instant. 

We also notice that we have $2|\mathcal{E}| > N$, which implies that $\Am$ is not full row-rank. As we will see shortly, this implies that the dual problem to~\eqref{eq:consensus-problem} is not strongly-convex. This generally is a problem in time-varying optimization, since then the (dual) optimizer trajectory is not unique. We explore next how to tackle this issue with a dual-regularized version of ADMM.

\subsection{Dual-regularized Problem}\label{sec:dual-regularization}

We now examine more closely the consensus problem~\eqref{eq:consensus-problem}. First, we can write its dual problem as
\begin{equation}\label{eq:dual-problem}
	\w^*(t_k) = \argmin_\w \left\{ d^{f_k}(\w) + d^h(\w) \right\}
\end{equation}
where we have set $f(\x; t_k) := \sum_{i=1}^N f_i(\x_i; t_k)$ and $h(\y) = \iota_{\ker(\Im - \Pm)}(\y) $, and where we have:
$$
	d^{f_k}(\w) = f^\star(\Am^\top \w; t_k) \quad \text{and} \quad d^h(\w) = h^\star(-\w).
$$

The dual problem is strongly convex provided that the cost $f(\cdot; t_k)$ is $L$-smooth and that $\Am$ is full row-rank. However, as noticed, for distributed problems, $\Am$ is rank deficient, and thus the dual problem is only convex, see also \cite[Remark~5]{bastianello_asynchronous_2019}. This is an issue for time-varying algorithms\footnote{Notice that the distributed ADMM derived by applying the \emph{Douglas-Rachford splitting} to the dual has provable linear convergence, see \cite{shi_linear_2014}. However, we purposely use the Peaceman-Rachford splitting applied to the regularized dual problem, since its convergence rate is better than that of the Douglas-Rachford, which is key in time-varying scenarios, see also~\cite{bastianello_asynchronous_2019}.}.  

In order to have a strongly convex dual problem, we introduce now a dual regularization, by substituting the linear constraints~\eqref{eq:lin-constraint} with the following
\begin{equation}\label{eq:regularized-constraints}
	\Am \x - \y = - (\epsilon/2) \w
\end{equation}
with $\epsilon > 0$. The corresponding dual problem is then (see \ifarxiv Appendix~\ref{app:regularized-dual}\else \cite[Appendix~A]{arxiv_version}\fi)
\begin{equation}\label{eq:dual-reg-problem}
	\w^*(\epsilon; t_k) = \argmin_\w \left\{ d^{f_k}(\w;\epsilon) + d^h(\w) \right\}
\end{equation}
with
$
	d^{f_k}(\w;\epsilon) = \frac{\epsilon}{2} \norm{\w}^2 + f^\star(\Am^\top \w; t_k),
$
and where we have overloaded the notation with: $\w^*(t_k) = \w^*(0;t_k)$ and $d^{f_k}(\w) = d^{f_k}(\w;0)$. Problem~\eqref{eq:dual-reg-problem} has now a $\epsilon$-strongly convex cost function.

Before moving on to solve~\eqref{eq:dual-reg-problem}, we examine its properties with respect to the original~\eqref{eq:dual-problem}, under the following assumption.

\smallskip \begin{assumption}\label{as:finite-solution}
The solutions of the original and regularized dual problems, \eqref{eq:dual-problem} and \eqref{eq:dual-reg-problem}, respectively, are finite, for each time instance $t_k$: $\max\{\|\w^*(t_k)\|, \|\w^*(\epsilon; t_k)\|\}\leq C$. 
\end{assumption} \smallskip

The following Lemma bounds the distance between the regularized dual solution and the solution(s) of the original problem around a neighborhood of $\epsilon = 0$.

\smallskip
\begin{lemma}[Regularization error]\label{lem:regularization-error}
Let the cost $f(\cdot; t_k)$ be in $\mathcal{S}_{\mu,L}(\reals^{nN})$ uniformly in time. Let $\w^*(\epsilon; t_k)$ be the optimal solution of the regularized dual problem~\eqref{eq:dual-reg-problem}, and $\w^*(t_k)$ be a solution to the original dual problem~\eqref{eq:dual-problem}. Then there exist a  $\psi_0 > 0$ for which for all $\epsilon>0$ close enough to $0$, we have
$$
	\norm{\w^*(\epsilon; t_k) - \w^*(t_k)} \leq (1 + \psi_0 \epsilon) \norm{\w^*(\epsilon; t_k)}.
$$
\end{lemma}
\smallskip
\begin{proof}
See \ifarxiv Appendix~\ref{app:proof-distance-optima}\else \cite[Appendix~B]{arxiv_version}\fi.
\end{proof}
\smallskip

Notice that by Lemma~\ref{as:finite-solution}, it holds that:
$$
	\norm{\w^*(\epsilon; t_k) - \w^*(t_k)} \leq C (1 + \psi_0 \epsilon).
$$
The bound provided by Lemma~\ref{lem:regularization-error} is not tight, in the sense that when $\epsilon \to 0$, the norm $\norm{\w^*(\epsilon; t_k) - \w^*(t_k)}$ does not go to zero. Nonetheless, the bound does characterize the fact that the smaller $\epsilon$ is, the smaller the regularization error is, guaranteeing that it is indeed a bounded error.

\subsection{Distributed Dual-regularized ADMM}

We are now ready to solve~\eqref{eq:dual-reg-problem} for a fixed time instant $t_k$. One can use the Peaceman-Rachford splitting to the regularized dual problem to find the corresponding dual-regularized ADMM that solves the primal and dual problem (see \ifarxiv Appendix~\ref{app:regularized-dual}\else \cite[Appendix~A]{arxiv_version}\fi for the details).

First, we introduce the \emph{auxiliary variables} $\z := \{ \z_{ij}, \z_{ji} \}_{(i,j) \in \mathcal{E}} \in \reals^{2n|\mathcal{E}|}$, two for each edge, from which we can compute the value of the dual variable $\w_{ij}$ for the respective consensus constraints. In particular, starting from an initial guess $\z^0 \in \reals^{2n |\mathcal{E}|}$, the dual-regularized ADMM corresponds to the following recursion for $\ell \in \mathbb{N}$:
\begin{subequations}\label{eq:regularized-admm}
\begin{align}
	&\x_k^\ell = \argmin_{\x} \big\{ f(\x; t_k) + \frac{\rho\delta}{2} \norm{\Am \x - \z_k^\ell / \rho}^2 \big\} \label{eq:reg-admm-x} \\
	&\w_k^\ell = \delta \left( \z_k^\ell - \rho \Am \x_k^\ell \right) \label{eq:reg-admm-w} \\
	&\y_k^\ell = \argmin_{\y} \big\{ h(\y) + \frac{\rho}{2} \norm{- \y - (2\w_k^\ell - \z_k^\ell) / \rho} \big\} \\
	&\uv_k^\ell = 2\w_k^\ell - \z_k^\ell + \rho \y_k^\ell \\
	&\z_k^{\ell+1} = \z_k^\ell + 2(\uv_k^\ell - \w_k^\ell)
\end{align}
\end{subequations}
where we have set $\delta := 1 / (1 + \epsilon \rho)$ for simplicity, and the subscript $_k$ denotes dependence from the problem sampled at time $t_k$. Recursions~\eqref{eq:regularized-admm} generate the sequence of primal and auxiliary variables $\{\x_k^{\ell}, \z_k^\ell \}_{\ell \in \mathbb{N}}$. Recursions~\eqref{eq:regularized-admm} can be shown to be implementable in a distributed fashion due to the particular distributed structure of the problem. In particular, it is possible to derive the following local updates (see \ifarxiv Appendix~\ref{app:distributed-admm}\else \cite[Appendix~D]{arxiv_version}\fi), called distributed dual-regularized ADMM. 

\begin{subequations}\label{eq:distributed-admm}

\smallskip
{\bf \sf Distributed Dual-regularized ADMM \hrule }

\begin{enumerate}
\item[\bf 0.] Initialization: node $i$ set $\z_{ij,k}^{0} = \zero$ for all $j \in \mathcal{N}_i$
\item[\bf 1.] Local update: node $i$ computes
\begin{equation}
	\!\!\!\!\!\!\!\!\!\!\!\!\x_{i,k}^\ell = \argmin_{\x_i \in \reals^n} \Big\{ f_i(\x_i;t_k) + \frac{\rho \delta d_i}{2} \norm{\x_i}^2 - \delta \langle \x_i, \sum_{j \in \mathcal{N}_i} \z_{ij, k}^\ell \rangle \Big\} \label{eq:admm-x}
\end{equation}
\item[\bf 2.] Communication step: node $i$ sends to neighbor $j \in \mathcal{N}_i$ the local variables $\x_{i,k}^\ell$ and $\z_{ij,k}^\ell$
\item[\bf 3.] Auxiliary update: using the information received from node $j \in \mathcal{N}_i$, node $i$ computes
\begin{equation}
	\z_{ij,k}^{\ell+1} = (2\delta - 1) \z_{ji,k}^\ell+ 2 \delta \rho \x_{j,k}^\ell \label{eq:admm-z}.
\end{equation}
\hrule
\end{enumerate}
\end{subequations}
\smallskip

As a consequence of \cite[Lemma~A.3]{bastianello_primal_2020}, the dual-regularized ADMM~\eqref{eq:regularized-admm}, and its distributed version, converges Q-linearly to the solution of the dual-regularized problem as
$$
	\norm{\w_k^{\ell+1} - \w^*(\epsilon; t_k)} \leq \zeta(\ell;\epsilon) \norm{\w_k^\ell - \w^*(\epsilon; t_k)}
$$
where $\zeta(\ell; \epsilon) \in (0,1)$ is defined in~\eqref{eq:zeta}, and R-linearly to the solution of the primal dual-regularized problem as
$$
	\norm{\x_k^\ell - \x^*(\epsilon, t_k)} \leq (\norm{\Am} / \mu ) \norm{\w_k^\ell - \w^*(\epsilon; t_k)}.
$$

\section{Distributed Prediction-Correction ADMM}\label{sec:pc-admm}

The previous sections have introduced a dual regularization technique to guarantee strong convexity of the dual problem, and thus linear convergence of the ADMM. Moreover, we have formulated the time-varying distributed problem of interest as a sequence of static problems that can be solved with ADMM. In this section, we briefly review the prediction-correction framework analyzed in \cite{bastianello_primal_2020}, and apply it to solve the distributed problem~\eqref{eq:consensus-problem}.

The proposed prediction-correction scheme is characterized by the following two steps:

$\bullet$ \emph{Prediction}: at time $t_k$, each node approximates the as yet unobserved local cost $f_{i,k+1}(\x_i) := f_i(\x_i;t_{k+1})$ using past observations of the cost; let $\hat{f}_{i,k+1}(\x)$ be such approximation, then the network solves
	\begin{equation}\label{eq:prediction-problem}
	\begin{split}
		&\min_{\x \in \reals^{nN}} \sum_{i=1}^N \hat{f}_{i,k+1}(\x_i) + \iota_{\ker(\Im - \Pm)}(\y) \\
		&\text{s.t.} \ \Am \x - \y = 0
	\end{split}
	\end{equation}
	which yields the prediction $\hat{\x}_{k+1}^*$. In practice, it is possible to compute only an approximation of $\hat{\x}_{k+1}^*$, denoted by $\hat{\x}_{k+1}$, by applying $\Np$ steps of the dual-regularized ADMM.
	
$\bullet$ \emph{Correction}: when, at time $t_{k+1}$, the nodes can observe $f_{i,k+1}(\x)$, they can correct the prediction computed at the previous step by solving:
	\begin{equation}\label{eq:correction-problem}
		\begin{split}
			&\min_{\x \in \reals^{nN}} \sum_{i=1}^N f_{i,k+1}(\x_i) + \iota_{\ker(\Im - \Pm)}(\y) \\
			&\text{s.t.} \ \Am \x - \y = 0
		\end{split}
		\end{equation}
	with initial condition equal to $\hat{\x}_{k+1}$. We will denote by $\x_{k+1}$ the (possibly approximate) correction computed by applying $\Nc$ of the dual-regularized ADMM.
\smallskip

During the prediction step, each node needs to approximate the future cost using the information \emph{locally} available up to time $t_k$. In particular, we employ the following Taylor-based expansion of $\nabla_{\x_i} f_{i,k+1}$ around $(\x_{i,k},t_k)$:
\begin{equation}
\begin{split}
	&\nabla_{\x_i} \hat{f}_{i,k+1}(\x_i) := \nabla_{\x_i} f_{i,k}(\x_{i,k}) \\ &\quad + \nabla_{\x_i \x_i} f_{i,k}(\x_{i,k}) (\x_i - \x_{i,k}) + \Ts \nabla_{t \x_i} f_{i,k}(\x_{i,k}).
\end{split}
\end{equation}
Notice that $\hat{f}_{i,k+1}$ is a quadratic function that inherits the same strong convexity and smoothness properties of $f_{i,k}$, since they have the same Hessian.

Algorithm~\ref{alg:prediction-correction} reports the pseudo-code for the proposed prediction-correction ADMM.

\begin{algorithm}[!ht]
\caption{Prediction-correction dual-regularized ADMM.}
\label{alg:prediction-correction}
\footnotesize
\begin{algorithmic}[1]
	\Require $\x_{i,0}$, horizons $\Np$ and $\Nc$, parameters $\epsilon$ and $\rho$.
	\For{$k=0,1,\ldots$}	
	\Statex\hspace{\algorithmicindent}{\color{blue}// time $t_k$ (prediction)}
	\For{$i \in \mathcal{V}$}
		\State Compute the local prediction function $\hat{f}_{i,k+1}(\x_i)$
	\EndFor
	\State Apply $\Np$ steps of the distributed dual-regularized ADMM (from $\ell = 1$ to $\ell = \Np$) to problem~\eqref{eq:prediction-problem}; that is, the nodes perform the steps 1)-3) $\Np$ times, choosing ${\z}_{ij,k+1}^0 = \z_{ij,k}$ and outputting $\hat{\z}_{ij,k+1} = {\z}_{ij,k+1}^{\Np}$ and $\hat{\x}_{i,k+1} = {\x}_{i,k+1}^{\Np}$; 
	\Statex\hspace{\algorithmicindent}{\color{blue}// time $t_{k+1}$ (correction)}
	\For{$i \in \mathcal{V}$}
		\State Observe the local cost function $f_{i,k+1}(\x_i)$
	\EndFor
	\State Apply $\Nc$ steps of the distributed dual-regularized ADMM (from $\ell = 1$ to $\ell = \Nc$) to problem~\eqref{eq:correction-problem}; that is, the nodes perform the steps 1)-3) $\Nc$ times, choosing $\z_{ij,k+1}^0 = \hat{\z}_{ij,k+1}$ and outputting $\z_{ij,k+1} = \z_{ij,k+1}^{\Nc}$ and $\x_{i,k+1} = \x_{i,k+1}^{\Nc}$; 
	\State Set $\x_{k+1}$ equal to the last iterate of ADMM
	\EndFor
\end{algorithmic}
\end{algorithm}

\subsection{Convergence analysis}

We now analyze the convergence properties of Algorithm~\ref{alg:prediction-correction} to find and track the solution trajectory $\x^*(t_k)$ of the original problem~\eqref{eq:consensus-problem}. We will use the following standard (in time-varying optimization) assumption. 

\smallskip \begin{assumption}\label{as:problem-properties}
\emph{(i)} The local costs $f_i : \reals^n \times \reals_+ \to \reals$ belong to $\mathcal{S}_{\mu,L}(\reals^n)$ uniformly in $t$. \emph{(ii)} There exists $C_0$ such that $\norm{\nabla_{t \x_i} f_i(\x_i;t_k)} \leq C_0$ for any $\x_i \in \reals^n$, $t \in \reals_+$. \emph{(iii)} The solution to~\eqref{eq:consensus-problem} is finite for any $k \in \mathbb{N}$.
\end{assumption} \smallskip

The following Lemma lists the properties of the regularized dual problem given that Assumption~\ref{as:problem-properties} holds for the primal.

\smallskip
\begin{lemma}\label{lemma-dual-prop}
Let Assumption~\ref{as:problem-properties} hold. The dual function $d^{f_k}(\w;\epsilon)$ is $\bar{\mu}(\epsilon) := \epsilon$-strongly convex and $\bar{L}(\epsilon) := \epsilon + d_\mathrm{M} / \mu$-smooth. Moreover, for any $\w$, $t_k$, $\epsilon$ it holds $\norm{\nabla_{t\w} d^{f_k}(\w,\epsilon)} \leq \sqrt{d_\mathrm{M}} C_0 / \mu =: \bar{C}_0$.
\end{lemma}
\smallskip
\begin{proof}
By \cite[Lemma~A.1]{bastianello_primal_2020} the dual function $f^{\star}(\Am^\top \w; t_k)$ is $\lambda_\mathrm{M}(\Am \Am^\top) / \mu$-smooth, and $\lambda_\mathrm{m}(\Am \Am^\top) / L$-strongly convex. But since $\lambda_\mathrm{m}(\Am \Am^\top) = 0$, then it is only convex. Therefore, adding the regularization term implies $\epsilon$-strong convexity and $ \epsilon + d_\mathrm{M} / \mu$-smoothness of $d^{f_k}$.

By \cite[Lemma~4.3]{bastianello_primal_2020} we have that $\norm{\nabla_{t\w} d^{f_k}(\w;\epsilon)} \leq \norm{\Am} C_0 / \mu$. Moreover, we have that $\Am \Am^\top = \operatorname{blk\,diag}(\one_{d_i \times d_i}) \otimes \Im_n$, and so the $i$-th diagonal block has eigenvalues $d_i$ and $0$. Since $\norm{\Am} = \sqrt{\lambda_\mathrm{M}(\Am \Am^\top)} = \sqrt{\max_i d_i}$, it follows $\norm{\Am} = \sqrt{d_\mathrm{M}}$.
\end{proof}

\smallskip

Notice that as a consequence of Lemma~\ref{lemma-dual-prop}, the condition number of $d^{f_k}$ is $\bar{\kappa}(\epsilon) := 1 + d_\mathrm{M} / (\mu \epsilon)$ independently of the time instant $t_k$. 

As mentioned above, the ADMM corresponds to the Peaceman-Rachford splitting (PRS) applied to the regularized dual. Since the dual is strongly convex and smooth, we know that the PR operator is Lipschitz continuous with constant \cite{giselsson_linear_2017}
$$
	\lambda(\epsilon) := \max\left\{ \Big\lvert \frac{1 - \rho \bar{L}(\epsilon)}{1 + \rho \bar{L}(\epsilon)} \Big\rvert, \Big\lvert \frac{1 - \rho \bar{\mu}(\epsilon)}{1 + \rho \bar{\mu}(\epsilon)} \Big\rvert \right\} \in (0,1).
$$

Moreover, given the fixed point $\z^*(\epsilon;t_k)$ of the PR operator, that is, $\z^*(\epsilon;t_k) = \refl_{\rho d^h} (\refl_{\rho d^f}(\z^*(\epsilon;t_k)))$, we can derive the solution to the dual with $\w^*(\epsilon;t_k) = \prox_{\rho d^{f_k}}(\z^*(\epsilon;t_k))$, and it holds 
$$
	\norm{\w - \w^*(\epsilon;t_k)} \leq \omega(\epsilon) \norm{\z - \z^*(\epsilon;t_k)},
$$
where $\z \in \reals^{2n|\mathcal{E}|}, \w = \prox_{\rho d^{f_k}}(\z)$ and with $\omega(\epsilon) := (1 + \rho \bar{L}(\epsilon)) / (1 + \rho \bar{\mu}(\epsilon))$.

Finally, we introduce the following notation that will be useful for the convergence results:
\begin{align}
	\zeta(\ell;\epsilon) &:= \left\{\begin{array}{lr} 1, & \textrm{for } \ell = 0, \\ \omega(\epsilon) \lambda(\epsilon)^\ell, & \textrm{otherwise} \end{array}  \right. \label{eq:zeta} \\
	\xi(\ell;\epsilon) &:= \left\{\begin{array}{lr} 0, & \textrm{for } \ell = 0, \\ 1 + \omega(\epsilon) \lambda(\epsilon)^\ell, & \textrm{otherwise} \end{array}  \right. \nonumber
\end{align}

The following result characterizes the convergence of the sequence of $\{ \w_k \}_{k \in \mathbb{N}}$ generated by the prediction-correction ADMM of Algorithm~\ref{alg:prediction-correction} to a neighborhood of the regularized dual optimal trajectory $\{ \w^*(\epsilon; t_k) \}_{k \in \mathbb{N}}$.

\begin{remark}
As one can notice, Algorithm~\ref{alg:prediction-correction} technically does not generate dual variables, since not necessary, yet one could generate those by implementing~\eqref{eq:reg-admm-w} after Step 1 in a distributed way. In particular, node $i$ can compute the dual variables $\w_{ij,k}$, $j \in \mathcal{N}_i$ as $\w_{ij,k} = \delta (\z_{ij,k} - \rho \x_{i,k})$, which only requires the local information $\z_{ij,k}$ and $\x_{i,k}$, without the need for additional communications.
\end{remark}

\smallskip
\begin{corollary}[Convergence to $\w^*(\epsilon; t_k)$]\label{cor:convergence-to-regularized-sol}
Let As.~\ref{as:problem-properties} hold. Let $\{ \w_k \}_{k \in \mathbb{N}}$ be the sequence of dual variables generated by the prediction-correction ADMM of Algorithm~\ref{alg:prediction-correction}. Choose $\Np$, $\Nc$ such that
\begin{equation}\label{eq:convergence-condition}
	\eta_1(\epsilon) := \zeta(\Nc;\epsilon) \left[ \zeta(\Np;\epsilon) + 2 \bar{\kappa}(\epsilon) \xi(\Np;\epsilon) \right] < 1.
\end{equation}
Then the trajectory $\{ \w_k \}_{k \in \mathbb{N}}$ converges to a neighborhood of $\{ \w^*(\epsilon; t_k) \}_{k \in \mathbb{N}}$ with radius upper bounded by
$$
	\limsup_{k \to \infty} \norm{\w_k - \w^*(\epsilon; t_k)} = \frac{\eta_0(\epsilon)}{1 - \eta_1(\epsilon)},
$$
where
$$
	\eta_0(\epsilon) := \zeta(\Nc;\epsilon) \frac{\bar{C}_0 \Ts}{\bar{\mu}} \left[ \zeta(\Np;\epsilon) + 2(1 + \bar{\kappa}(\epsilon) \xi(\Np;\epsilon)) \right].
$$
\end{corollary}
\smallskip
\begin{proof}
This is a consequence of \cite[Theorem~3.10]{bastianello_primal_2020} applied to the regularized dual problem.
\end{proof}
\smallskip

The following result characterizes the convergence of $\{ \w_k \}_{k \in \mathbb{N}}$ in terms of the original optimal trajectory.

\smallskip
\begin{proposition}[Convergence to $\w^*(t_k)$]\label{pr:dual-convergence}
Let As.~\ref{as:finite-solution}-\ref{as:problem-properties} hold. Let $\{ \w_k \}_{k \in \mathbb{N}}$ be the sequence of dual variables generated by the prediction-correction ADMM of Algorithm~\ref{alg:prediction-correction}, and choose $\Np$, $\Nc$ such that~\eqref{eq:convergence-condition} holds. Then, for a small enough $\epsilon>0$, $\{ \w_k \}_{k \in \mathbb{N}}$ converges to a neighborhood of the original solution $\{ \w^*(t_k) \}_{k \in \mathbb{N}}$ with radius upper bounded by
$$
	\limsup_{k \to \infty} \norm{\w_k - \w^*(t_k)} = \frac{\eta_0(\epsilon)}{1 - \eta_1(\epsilon)} + C (1 + \psi_0 \epsilon).
$$
\end{proposition}
\smallskip
\begin{proof}
The result follows by the triangle inequality:
$$
	\norm{\w_k \!-\! \w^*(t_k)} \!\leq\! \norm{\w_k \!-\! \w^*(\epsilon; t_k)} \!+\! \norm{\w^*(t_k) \!-\! \w^*(\epsilon; t_k)},
$$
from the bound in Lemma~\ref{lem:regularization-error} to the regularization error, and taking the limit and using Corollary~\ref{cor:convergence-to-regularized-sol}.
\end{proof}
\smallskip

Finally, we can characterize the convergence of the prediction-correction ADMM in terms of the primal variable.

\begin{corollary}[Convergence to $\x^*(t_k)$]\label{cor:primal-convergence}
Let As.~\ref{as:finite-solution}-\ref{as:problem-properties} hold. Let $\{ \x_k \}_{k \in \mathbb{N}}$ be the sequence of primal variables generated by the prediction-correction ADMM of Algorithm~\ref{alg:prediction-correction}, and choose $\Np$, $\Nc$ such that~\eqref{eq:convergence-condition} holds. Then, for small enough $\epsilon>0$, the primal trajectory converges to a neighborhood of the original optimal trajectory $\{ \x^*(t_k) \}_{k \in \mathbb{N}}$ with radius upper bounded as:
$$
	\limsup_{k \to \infty} \norm{\x_k - \x^*(t_k)} = \frac{\norm{\Am}}{\mu} \left[ \frac{\eta_0(\epsilon)}{1 - \eta_1(\epsilon)} + C (1 + \psi_0 \epsilon) \right].
$$
\end{corollary}
\smallskip
\begin{proof}
The result follows by combining Proposition~\ref{pr:dual-convergence} with the following bound
$$
	\norm{\x_k - \x^*(t_k)} \leq (\norm{\Am} / \mu) \norm{\w_k - \w^*(t_k)},
$$
which can be derived for the dual-regularized ADMM along the lines of \cite[Lemma~A.3]{bastianello_primal_2020}.
In particular, by the strong convexity of $f(\cdot; t_k)$, we know that
\begin{multline}\label{eq:str-cvx}
	\mu \norm{\x - \x^*(t_k)}^2 \leq \\ \langle \nabla_\x f(\x; t_k) - \nabla_\x f(\x^*(t_k); t_k), \x - \x^*(t_k)\rangle
\end{multline}
for any $\x \in \reals^{nN}$. Moreover, by the KKT conditions of the original problem~\eqref{eq:consensus-problem}, it must hold $\nabla_\x f(\x^*(t_k); t_k) = \Am^\top \w^*(t_k)$.

Imposing the first-order optimality condition for~\eqref{eq:reg-admm-x} at time $\ell = \Nc$ we have
$$
	\nabla_\x f(\x_k; t_k) - \delta \Am^\top \left( \z_k - \rho \Am \x_k \right) = \zero
$$
and, recalling by~\eqref{eq:reg-admm-w} that $\w_k = \delta \left( \z_k - \rho \Am \x_k \right)$, yields $\nabla_\x f(\x_k; t_k) = \Am^\top \w_k$.

Subtracting $\nabla_\x f_k(\x_k) = \Am^\top \w_k$ and $\nabla_\x f(\x^*(t_k); t_k) = \Am^\top \w^*(t_k)$, and substituting them into~\eqref{eq:str-cvx} we get
\begin{align*}
	\mu \norm{\x_k - \x^*(t_k)}^2 &\leq \langle \Am^\top (\w_k - \w^*(t_k)), \x_k - \x^*(t_k) \rangle \\
	&\leq \norm{\Am} \norm{\w_k - \w^*(t_k)} \norm{\x_k - \x^*(t_k)}
\end{align*}
which proves the thesis.
\end{proof}
\smallskip

\subsection{Trade-offs}
It is interesting now to illustrate the trade-off, mediated by $\epsilon$, between convergence rate and the asymptotic error achieved by the proposed algorithm. On one hand, low values of $\epsilon$ imply algorithms closer to the original problem, on the other hand, high values of $\epsilon$ are better for a more favorable convergence rate. This is translated into a larger or smaller radius for the asymptotical error.

For simplicity, in the following we take $\rho = 1 / \sqrt{\bar{\mu}(\epsilon) \bar{L}(\epsilon)}$, which maximizes $\lambda(\epsilon)$. As a consequence we have:
$$
	\omega(\epsilon) = \sqrt{\bar{\kappa}(\epsilon)}, \quad \lambda(\epsilon) = \Bigg\lvert \frac{1 - \sqrt{\bar{\kappa}(\epsilon)}}{1 + \sqrt{\bar{\kappa}(\epsilon)}} \Bigg\rvert,
$$
and using these facts we get (assuming that $\ell > 0$):
$$
	\zeta(\ell;\epsilon) = \sqrt{\bar{\kappa}(\epsilon)} \Bigg\lvert \frac{1 - \sqrt{\bar{\kappa}(\epsilon)}}{1 + \sqrt{\bar{\kappa}(\epsilon)}} \Bigg\rvert^\ell, \qquad 
	\xi(\ell;\epsilon) = 1 + \zeta(\ell;\epsilon).
$$
Since $\sqrt{\bar{\kappa}(\epsilon)}$ is monotonically decreasing as $\epsilon$ increases, and $\lim_{\epsilon \to +\infty} \sqrt{\bar{\kappa}(\epsilon)} = 1$, then it follows that $\zeta(\ell;\epsilon)$ and $\xi(\ell;\epsilon)$ decrease as well when $\epsilon$ grows. In turn, this implies that $\eta_0$ and $\eta_1$ decrease as well as $\epsilon$ increases.

Overall, the asymptotic error term $\eta_0(\epsilon) / (1 - \eta_1(\epsilon))$ due to the prediction-correction scheme grows smaller as larger values of the regularization constant $\epsilon$ are chosen. On the other hand, however, the regularization error itself increases when $\epsilon$ does so.
Therefore we can observe that $\epsilon$ mediates a trade-off between the two terms in the asymptotic error.

\subsection{Communication complexity}
Recalling the dual-regularized ADMM of~\eqref{eq:admm-x}-\eqref{eq:admm-z}, we can observe that node $i$ sends at each iteration $d_i$ packets, one to each of its neighbors. As a consequence, in Algorithm~\eqref{alg:prediction-correction} at each sampling time $t_k$ the total number of communications performed by node $i$ is equal to $(\Np + \Nc) d_i$.

This highlights a further trade-off, between the communication complexity and tracking error, which respectively increase and decrease as $\Np$ and $\Nc$ grow larger.

\section{Numerical results}
In this section, we present numerical results showcasing the performance of the proposed algorithm on the distributed optimization problem characterized by the local costs \cite{dallanese_optimization_2019}:
$$
	f_i(x_i;t) = \frac{1}{2} \norm{x_i - b_i(t)}^2 + \log\left( 1 + \exp(x_i - a_i) \right)
$$
where $b_i(t) = A \cos((\nu - 1) t + \varphi_i)$, $a_i \sim \mathcal{U}[-10,10]$, $\varphi_i \sim \mathcal{U}[0,2\pi)$, $A = 2.5$, $\nu = \pi / 80$. The network is a random geometric graph with $N = 25$ nodes. The simulations were implemented using \texttt{tvopt} \cite{tvopt}.

We compare the proposed dual-regularized PC-ADMM with \emph{(i)} the prediction-correction \emph{dual decomposition} of \cite{simonetto_dual_2019} applied to $\min_\x f(\x; t_k)$ s.t. $(\Im - \Wm) \x = \zero$, with $\Wm$ a doubly stochastic matrix designed by the Metropolis-Hastings rule, and \emph{(ii)} the prediction-correction \emph{gradient method} \cite{bastianello_primal_2020}, applied to $\min_\x f(\x; t_k) + (1/2\alpha) \x^\top (\Im - \Wm) \x$, with $\alpha > 0$ a suitable step-size. Notice that the dual decomposition algorithm is imposing \emph{exact consensus constraints} in the static case, while the other two methods approximate them with a primal or dual regularization.
The parameters of the three algorithms were hand-tuned to achieve the best asymptotic error, and in particular for the dual regularized PC-ADMM we chose $\epsilon = 10^{-3}$ and $\rho = 1.06 \times 10^4$.

\begin{figure}[!ht]
\centering
\includegraphics[width=0.35\textwidth]{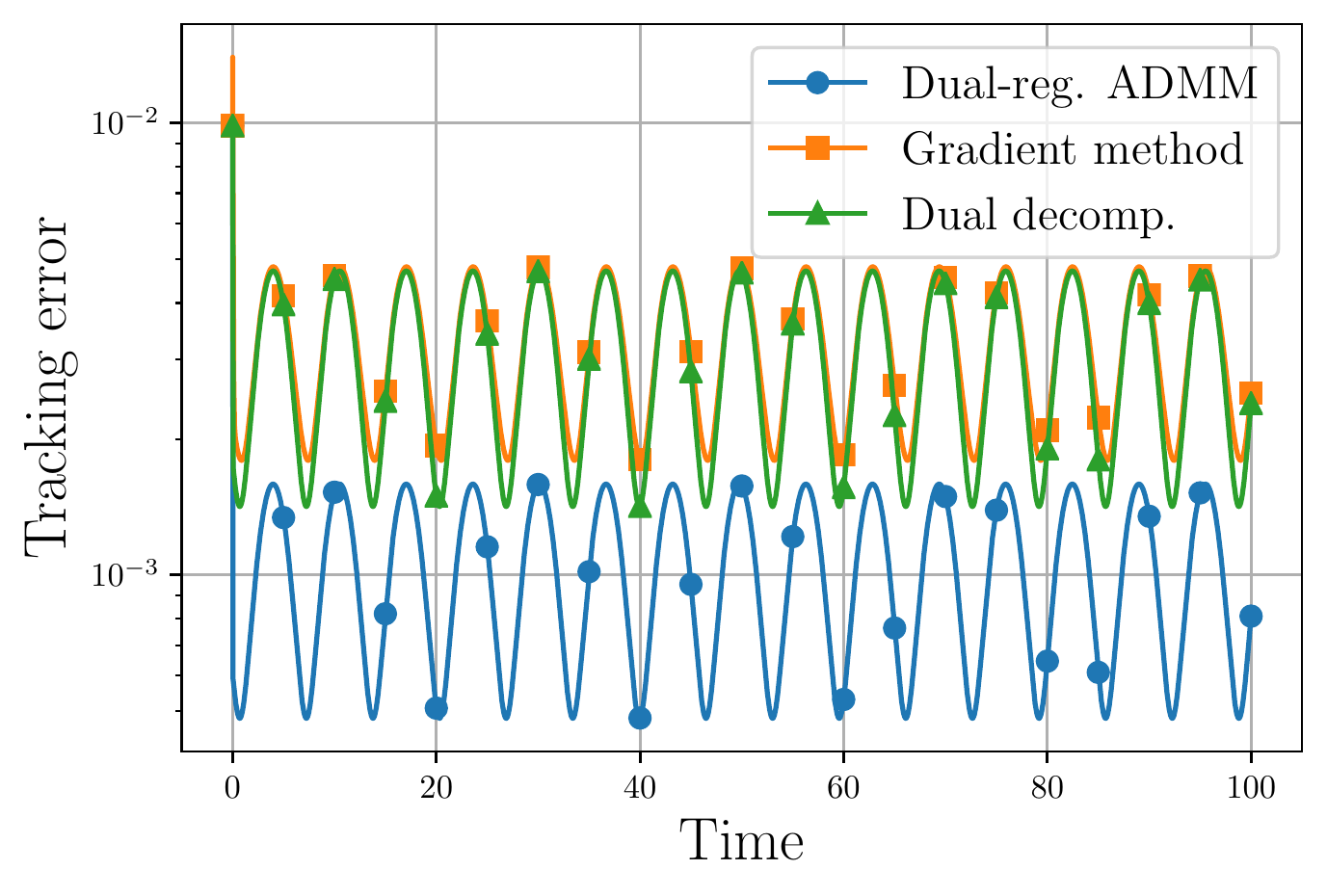}
\caption{Error trajectory comparison, with $\Np, \Nc = 5$.}
\label{fig:error}
\end{figure}

Figure~\ref{fig:error} depicts the trajectory of the error computed as $\norm{\x_k - \x^*(t_k)} / N$, where $\x^*(t_k)$ is the primal solution to the non-regularized problem.

\begin{figure}[!ht]
\centering
\includegraphics[width=0.35\textwidth]{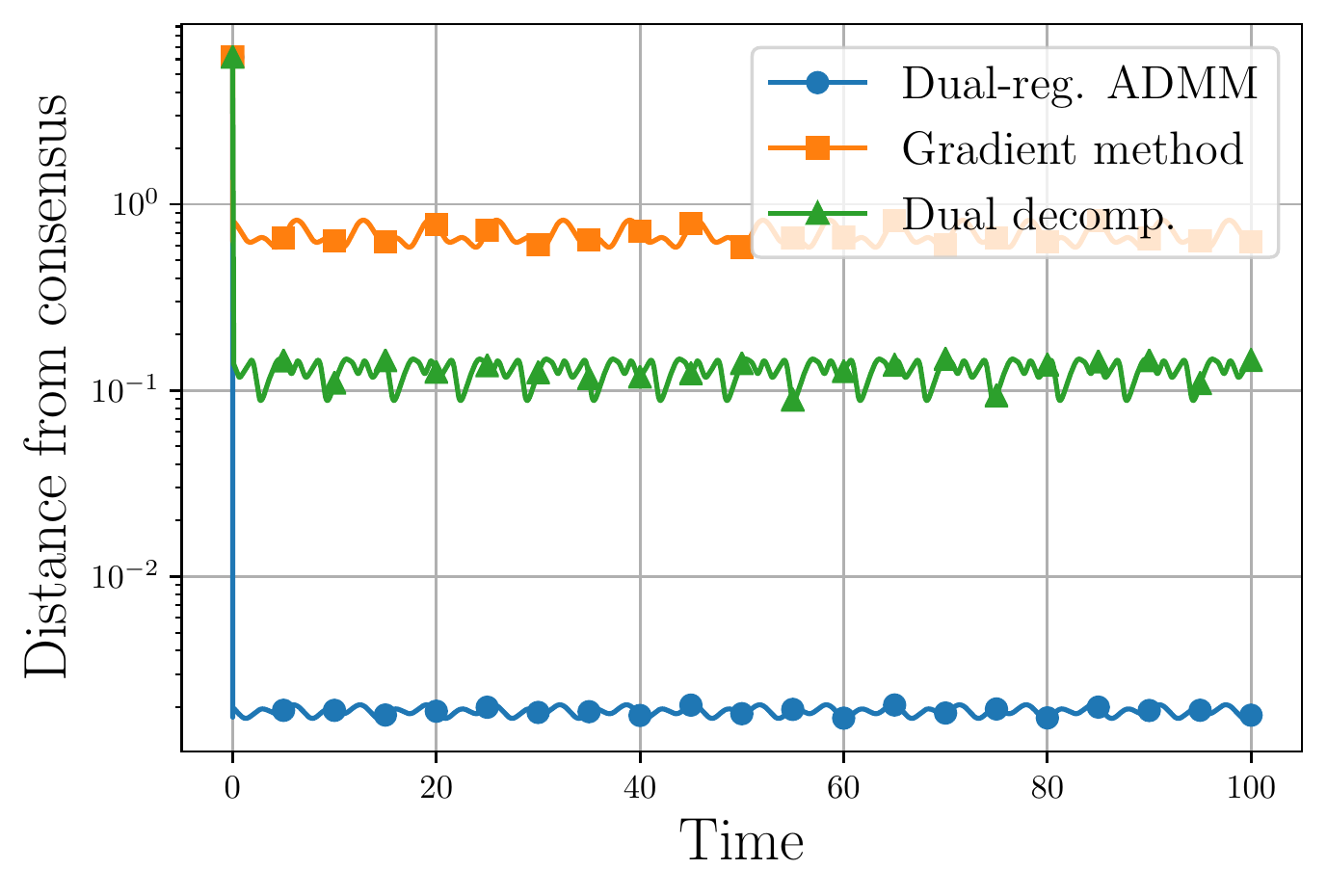}
\caption{Distance from consensus comparison, with $\Np, \Nc = 5$.}
\label{fig:consensus}
\end{figure}

Figure~\ref{fig:consensus} compares the evolution of the distance from consensus, computed as $\norm{\x_k - \one (\one^\top / N) \x_k}$.

As we can see from both results, the dual-regularized ADMM achieves the best performance in terms of asymptotic error and distance from consensus, even though the consensus constraints are not enforced exactly. Moreover, since the proposed algorithm outperforms the gradient method, it appears that a dual regularization is better than the primal regularization $(1/2\alpha) \x^\top (\Im - \Wm) \x$. This example further demonstrates how adding regularization terms in time-varying setting does not necessarily incurs in an accuracy trade-off.

\section{Conclusions}

In this paper, we have proposed a dual-regularized, prediction-correction ADMM to solve time-varying distributed optimization problems. On the one hand, the dual regularization ensures that the convergence of the algorithm is linear, and on the other, the prediction-correction scheme allows to efficiently track the optimal trajectory of the problem, up to a bounded error.


\ifarxiv

\appendix

\subsection{Regularized dual problem}\label{app:regularized-dual}
Similarly to \cite[Appendix~A.2]{peng_arock_2016} for the non-regularized ADMM, we define the Lagrangian of problem~\eqref{eq:consensus-problem} with the modified constraints~\eqref{eq:regularized-constraints}:
\begin{equation*}
	\mathcal{L}_k(\x,\y;\w) := f(\x; t_k) + h(\y) - \langle \w, \Am \x - \y  + (\epsilon/2) \w \rangle
\end{equation*}
and we compute the dual problem by minimizing $\mathcal{L}_k(\x,\y;\w)$ in $\x$ and $\y$. We have
\begin{multline*}
	\min_{\x, \y}\,\, \mathcal{L}_k(\x,\y;\w) = \min_\x \left\{ f(\x; t_k) - \langle \w, \Am \x \rangle \right\} \\ + \min_\y \left\{ h(\y) + \langle \w, \y \rangle \right\}  + \frac{\epsilon}{2} \norm{\w}^2 \\
	= f^\star(\Am^\top \w; t_k) + \frac{\epsilon}{2} \norm{\w}^2 + h^\star(- \w)
\end{multline*}
where we used the definition of convex conjugate. \endproof

\subsection{Proof of Lemma~\ref{lem:regularization-error}}\label{app:proof-distance-optima}

Let $\w_k^*(\epsilon) := \w^*(\epsilon; t_k)$ for convenience. First of all, given that $f$ is strongly convex, then the dual function is differentiable~\cite[Proposition~4.2]{bastianello_primal_2020}, and specifically,
\begin{equation}\label{eq.diff.dual}
\nabla_\w d^{f_k}(\w;\epsilon) = \epsilon \w + \nabla_\w f_k^\star(\Am^\top \w).
\end{equation}
Second, we notice that the dual optimal solution $\w_k^*(\epsilon)$ can be interpreted as the solution mapping $\epsilon \mapsto \w_k^*(\epsilon)$ of the generalized equation
\begin{equation}\label{gen}
	\nabla_\w d^{f_k}(\w;\epsilon) + \partial d^h(\w) \ni \zero,
\end{equation}
parametrized by $\epsilon$. The angle of attack for the proof is then bounding the Lipschitz continuity constant of this solution mapping using \cite[Theorem~2B.5]{dontchev2014implicit}.
We verify that the assumptions of \cite[Theorem~2B.5]{dontchev2014implicit} hold in this scenario, and the spell out the consequences of the Theorem.

Assumption (a) requires that $\nabla_\w d^{f_k}(\w;\epsilon)$ is continuous in $\epsilon$, and that the linearization, 
\begin{multline*}
\nabla_\w \hat{d}^{f_k}(\w) = \nabla_\w d^{f_k}(\w_k^*(\epsilon);\epsilon) +\\ \nabla_{\w\w} d^{f_k}(\w_k^*(\epsilon);\epsilon) (\w - \w_k^*(\epsilon))
\end{multline*}
is a strict estimator of $\nabla_\w d^{f_k}(\w;\epsilon)$ around $(\w_k^*(\epsilon);\epsilon)$ with constant $\mu$. Both are true, and in particular the linearization has $\mu=0$ (note in addition that $\nabla_{\w\w} d^{f_k}(\w_k^*(\epsilon);\epsilon)$ is well defined everywhere~\cite[Prop. 4.2]{bastianello_primal_2020}).

Assumption (b) requires that the inverse mapping $G_k^{-1}$ of $ G_k(\w):= \nabla_\w \hat{d}^{f_k}(\w) + \partial d^h(\w)$, for which $G_k(\w_k^*(\epsilon)) \ni \zero$, is Lipschitz continuous around $0$. But by definition of inverse mapping we have that
$$
	G_k^{-1}(\z) = \left\{ \w \ | \ \z \in \nabla_\w \hat{d}^{f_k}(\w) + \partial d^h(\w) \right\}
$$
which means that $G_k^{-1}(\z)$ is the solution mapping of the generalized equation $\z \in \nabla_\w \hat{d}^{f_k}(\w) + \partial d^h(\w)$ parametrized by $\z$. Since $\hat{d}^{f_k}$ is $\epsilon$-strongly convex, \cite[Theorem~D.1]{bastianello_primal_2020} implies that indeed $G_k^{-1}$ is $\epsilon^{-1}$-Lipschitz continuous everywhere.

Therefore, by \cite[Theorem~2B.5]{dontchev2014implicit}, the mapping~\eqref{gen} is single-valued and locally Lipschitz continuous around $\epsilon$ as
\begin{multline*}
\norm{\w_k^*(\epsilon) - \w_k^*(\epsilon')} \leq  (\epsilon^{-1} + \psi) \times \\ \norm{\nabla_\w d^{f_k}(\w_k^*(\epsilon);\epsilon) - \nabla_\w d^{f_k}(\w_k^*(\epsilon);\epsilon')},
\end{multline*}
where the locality is measured by $\psi > 0$. In particular, one can always choose a $\psi$, say $\psi_0$, such that $\epsilon'=0$ and $\epsilon>0$ in the neighborhood of $0$, and therefore by~\eqref{eq.diff.dual},
\begin{equation*}
\norm{\w_k^*(\epsilon) - \w_k^*(0)} \leq (1 + \psi_0 \epsilon) \norm{\w_k^*(\epsilon)},
\end{equation*}
from which the thesis follows. \endproof

\subsection{Regularized ADMM}\label{app:regularized-admm}
Following the derivation in \cite[Appendix~A]{davis_convergence_2016}, we now show that the Peaceman-Rachford splitting applied to the regularized dual problem is equivalent to~\eqref{eq:regularized-admm}.

The PRS is described by the updates, with $\ell \in \mathbb{N}$
\begin{subequations}
\begin{align}
	\w_k^\ell &= \prox_{\rho d^{f_k}}(\z_k^\ell) \label{eq:prs-w}\\
	\uv_k^\ell &= \prox_{\rho d^h}(2\w_k^\ell - \z_k^\ell) \label{eq:prs-u} \\
	\z_k^{\ell+1} &= \z_k^\ell + 2(\uv_k^\ell - \w_k^\ell)
\end{align}
\end{subequations}
and the aim is to show that~\eqref{eq:prs-w} is equivalent to
\begin{subequations}
\begin{align}
	&\x_k^\ell = \argmin_{\x} \left\{ f(\x; t_k) + \frac{\rho\delta}{2} \norm{\Am \x  - \z_k^\ell / \rho}^2 \right\} \\
	&\w_k^\ell = \delta \left( \z_k^\ell - \rho \Am \x_k^\ell \right);
\end{align}
\end{subequations}
the same derivation holds for~\eqref{eq:prs-u}.

By definition of the dual function $d^{f_k}$ and of proximal operator, to compute~\eqref{eq:prs-w} we need to solve the following minimization
\begin{align*}
	&\min_\w \left\{ d^{f_k}(\w) + \frac{1}{2\rho} \norm{\w - \z_k^\ell}^2 \right\} = \\
	&= \min_\w \left\{ f^\star(\Am^\top \w; t_k) + \frac{\epsilon}{2} \norm{\w}^2 + \frac{1}{2\rho} \norm{\w - \z_k^\ell}^2 \right\} \\
	&= \min_\w \max_\x \Big\{\! \langle \w, \Am \x \rangle\! - \!f(\x; t_k) \!+\! \frac{\epsilon}{2} \norm{\w}^2 \!+\! \frac{1}{2\rho} \norm{\w \!-\! \z_k^\ell}^2 \!\Big\} \\
	&= \max_\x \min_\w \Big\{\! \langle \w, \Am \x \rangle\! - \!f(\x; t_k) \!+\! \frac{\epsilon}{2} \norm{\w}^2 \!+\! \frac{1}{2\rho} \norm{\w \!-\! \z_k^\ell}^2 \!\Big\}
\end{align*}
where we used the definition of convex conjugate to the derive the second to last equality. Imposing the first order optimality condition for the innermost minimization yields
\begin{equation}\label{eq:w-update}
	\w_k^\ell = \delta \left( \z_k^\ell - \rho \Am \x  \right).
\end{equation}
Substituting~\eqref{eq:w-update} into the minimization problem and rearranging the terms gives
\begin{align*}
	\max_\x& \Big\{\! - \!f(\x; t_k) \!-\! \frac{\rho\delta}{2} \norm{\Am \x }^2 \!+\! \delta \langle \z_k^\ell , \Am \x \rangle \!+\! \frac{1-\delta}{2\rho} \norm{\z_k^\ell}^2 \!\Big\} \\
	&= - \min_\x \left\{ f(\x; t_k) + \frac{\rho\delta}{2} \norm{\Am \x  - \z_k^\ell / \rho}^2 \right\}
\end{align*}
whose argument yields the desired update for $\x$. The same computations can be traced to compute the proximal operator of $d^h$. \endproof

\subsection{Distributed ADMM}\label{app:distributed-admm}
This derivation follows closely \cite[Appendix~C-A]{bastianello_asynchronous_2019}. By definition of $\Am$, it is possible to see that $\norm{\Am \x}^2 = \x^\top (\operatorname{diag}{d_i} \otimes \Im_n) \x$, and moreover that $[\langle \Am^\top \z_k^\ell, \x \rangle]_i = \langle \sum_{j \in \mathcal{N}_i} \z_{ij,k}^\ell, \x_i \rangle$. Therefore we can write
$$
	\x_k^\ell \!=\! \argmin_\x \sum_{i=1}^N \Big\{\! f_{i}(\x_i; t_k) + \frac{\rho\delta d_i}{2} \norm{\x_i}^2 \!- \!\delta \langle \sum_{j \in \mathcal{N}_i} \z_{ij,k}^\ell, \x_i \rangle\! \Big\}.
$$
Since the problem is separable, the update~\eqref{eq:admm-x} for $\x_{i,k}^\ell$ follows.
Using the fact that $h$ is the indicator function of $\ker(\Im - \Pm)$ it is possible to derive (see \cite[Appendix~C-A]{bastianello_asynchronous_2019}). Indeed, imposing the KKT conditions for
$$
	\y_k^\ell = \argmin_{\y \in \ker(\Im - \Pm)} \left\{ \frac{\rho}{2} \norm{- \y - (2\w_k^\ell - \z_k^\ell) / \rho}^2 \right\}
$$
yields
\begin{align*}
	\y_k^\ell = \frac{1}{\rho} \left( (\Im - \Pm) \pmb{\nu} - (2\w_k^\ell - \z_k^\ell) \right) \quad \text{and} \quad \y_k^\ell = \Pm \y_k^\ell
\end{align*}
where $\pmb{\nu}$ are the Lagrange multipliers. Substituting the first KKT condition into the right-hand side of the second KKT condition gives
$$
	\y_k^\ell = \frac{1}{\rho} \left( - (\Im - \Pm) \pmb{\nu} - \Pm (2\w_k^\ell - \z_k^\ell) \right)
$$
where we used the property of permutations matrices $\Pm^2 = \Im$. Summing the two equations yields
\begin{equation}\label{eq:y-update}
	\y_k^\ell = - \frac{1}{2\rho} (\Im + \Pm) (2\w_k^\ell - \z_k^\ell)
\end{equation}
and substituting~\eqref{eq:y-update} into the update for $\uv$ yields:
\begin{equation}\label{eq:u-update}
	\uv_k^\ell = 2\w_k^\ell - \z_k^\ell + \rho \y_k^\ell = \frac{1}{2} (\Im - \Pm) (2\w_k^\ell - \z_k^\ell).
\end{equation}

Finally, using~\eqref{eq:u-update} and $\w_k^\ell = \delta (\z_k^\ell - \rho \Am \x_k^\ell)$ into the update for $\z$ yields
\begin{equation}
	\z_k^{\ell+1} = \z_k^\ell + 2(\uv_k^\ell - \w_k^\ell) = (1 - 2\delta) \Pm \z^\ell + 2\rho \delta \Pm \Am \x^\ell
\end{equation}
and using the definition of $\Pm$ we get~\eqref{eq:admm-z}. \endproof

\addtolength{\textheight}{-16.6cm}

\fi

\bibliographystyle{IEEEtran}
\bibliography{IEEEabrv,references}

\end{document}